\documentclass[a4paper, 12pt, onecolumn]{article}
\usepackage{geometry}

\usepackage{amsmath}
\usepackage{amsthm}
\usepackage{amsfonts}
\usepackage{bbm}
\usepackage{CJK}
\usepackage{fancyhdr}
\usepackage{graphicx}
\usepackage{geometry}
\usepackage{psfrag}
\usepackage{amsfonts,amsmath,amsthm, amssymb}
\usepackage{latexsym, euscript, epic, eepic}
\usepackage{time}
\usepackage{txfonts}
\usepackage{colortbl}
\usepackage{stmaryrd}
\usepackage{mathrsfs}
\usepackage{txfonts}
\usepackage{amsfonts}
\usepackage{color}
\usepackage{lineno}
\usepackage[square, comma, sort&compress, numbers]{natbib}
\usepackage{indentfirst,latexsym,bm}

 \numberwithin{equation}{section}

\begin{document}
\newtheorem{theorem}{Theorem}[section]
\newtheorem{proposition}[theorem]{Proposition}
\newtheorem{remark}[theorem]{Remark}
\newtheorem{corollary}[theorem]{Corollary}
\newtheorem{definition}{Definition}[section]
\newtheorem{lemma}[theorem]{Lemma}
\newtheorem{conjecture}[theorem]{Conjecture}
\newcommand{\IN}{\mathbb{N}}
\newcommand{\IR}{\mathbb{R}}
\newcommand{\IC}{\mathbb{C}}
\newcommand{\e}{\mbox{e}}
\newcommand{\arth}{\rm arth}
\newcommand{\K}{\mathscr{K}}
\newcommand{\E}{\mathscr{E}}

\title{Sharp Approximations for the Ramanujan
Constant
\thanks{This research is supported by NSF of China (Grant No.11171307) and Zhejiang Provincial NSF of China (Grant No.LQ17A010010).}}
\author{\normalsize Song-Liang Qiu \footnote{Corresponding author. E-mail address: sl$\_$qiu@zstu.edu.cn.}, Xiao-Yan Ma and Ti-Ren Huang\\
\scriptsize (Department of Mathematics, Zhejiang Sci-Tech University, Hangzhou 310018, China)}
\date{}
\maketitle
\fontsize{12}{22}\selectfont\small

\paragraph{Abstract:}
 In this paper, the authors present sharp approximations in terms of sine function and
polynomials for the so-called Ramanujan constant (or the Ramanujan $R$-function) $R(a)$, by showing some monotonicity, concavity and convexity properties of certain combinations defined in terms of $R(a)$, $\sin(\pi a)$ and polynomials. Some properties of the Riemann zeta function and its related special sums are presented, too.
\\[10pt]
\\[10pt]

\section{\normalsize Introduction}\label{sec:bd}

For real numbers $x, y>0$, the gamma, beta and psi functions
are defined as
\begin{eqnarray}\label{Gam-Psi}
\Gamma(x)=\int_0^\infty t^{x-1}e^{-t}dt,
~B(x,y)=\frac{\Gamma(x)\Gamma(y)}{\Gamma(x+y)},
~\psi(x)=\frac{\Gamma'(x)}{\Gamma(x)},
\end{eqnarray}
respectively. (For their basic properties, cf. \cite{1,2,5,14}.) Let
$\gamma=0.5772156649\cdots$ be the Euler-Mascheroni constant, and for $a\in(0,1)$, let $R(a)$ be defined by
\begin{eqnarray}\label{R}
R(a)=-2\gamma-\psi(a)-\psi(1-a)
\end{eqnarray}
which is called the Ramanujan constant in literature although it is actually a function of $a$ and probably better to call $R(a)$ the Ramanujan $R$-function (cf. \cite{11}). By the symmetry, we may assume that $a\in(0,1/2]$ in (\ref{R}).
It is well known that $R(a)$ is essential in some fields of mathematics
such as the zero-balanced Gaussian hypergeometric functions
${}_2F_1(a,1-a;1;z)$, the theories of Ramanujan's modular equations and
quasiconformal mappings, and the properties of $R(a)$ are
indispensable for us to show the properties of
 ${}_2F_1(a,1-a;1;z)$ and the functions appearing in generalized Ramanujan's modular
equations. On the other hand, $R(a)$ and the function
\begin{eqnarray}\label{B}
B(a)\equiv B(a,1-a)=\Gamma(a)\Gamma(1-a)=\frac{\pi}{\sin(\pi a)} ~(0<a<1)
\end{eqnarray}
are often simultaneously appear in the study of
the properties and applications of $R(a)$, and we often need to compare $R(a)$ with
$B(a)$. In \citep[Section 1]{11}, such kind of importance and applications of $R(a)$,
and the relation between $R(a)$ and $B(a)$ were described in details. (See also \cite{2,4,5,6,7,8,9,10,12,15,16,17}.)

Some authors have obtained some properties, including lower and upper bounds, for $R(a)$. In \cite{11}, for instance, power series expansion, integral representation and bounds were obtained for the difference $R(a)-B(a)$.
Some related studies showed that $B(a)$ is one of good approximation functions for $R(a)$, and we often require the properties of certain combinations defined in terms of $R(a)$, $B(a)$, $a(1-a)$ and other polynomials. (See \cite{3,8,10,11,15,16,17}.)

The main purpose of this paper is to show some monotonicity, convexity and concavity
 properties of certain combinations defined in terms of $R(a)$, $B(a)$ and polynomials, by which sharp
approximations given by $B(a)$ and
polynomials are obtained for $R(a)$. In addition, we shall also show
some properties of the Riemann zeta function
\begin{eqnarray}\label{zeta}
\zeta(s)=\sum_{k=1}^{\infty}\frac{1}{k^s}, Re~s>1
\end{eqnarray}
and its related special sums \citep[23.2]{1}
\begin{eqnarray}\label{lam-beta}
\lambda(n+1)=\sum_{k=0}^{\infty}\frac{1}{(2k+1)^{n+1}},
~\eta(n)=\sum_{k=1}^{\infty}(-1)^{k-1}\frac{1}{k^n},
~\beta(n)=\sum_{k=0}^{\infty}(-1)^k\frac{1}{(2k+1)^n}, n\in\IN.
\end{eqnarray}

For the later use, we record the following identities and special values
\begin{eqnarray}\label{lambda(n+1)}
\lambda(n+1)&=&\left(1-2^{-n-1}\right)\zeta(n+1), ~
\eta(n)=\left(1-2^{1-n}\right)\zeta(n), \label{lambda(n+1)}\\
\zeta(2)-\frac{\pi^2}{6}&=&\zeta(4)-\frac{\pi^4}{90}=\lambda(2)-\frac{\pi^2}{8}
=\lambda(4)-\frac{\pi^4}{96}=\beta(1)-\frac{\pi}{4}=\beta(3)-\frac{\pi^3}{32}=0.\label{zeta(2)}
\end{eqnarray}
(See \citep[23.2.19-20, 23.2.24-25 \& 23.2.28-31]{1}).

Throughout this paper, the ``zero-order derivative"
$\varphi^{(0)}$ will be understood as $\varphi$ itself for real one-variable function $\varphi$, and we always let
\begin{eqnarray}
&{}&
\begin{cases}
a_0=1, a_1=-2\eta(2)-1=-1-\pi^2/6=-2.644934\cdots, a_2=2\zeta(3)=2.404113\cdots,\\
a_n=\left[1+(-1)^n \right][\zeta(n+1)-\zeta(n-1)]+\left[1-(-1)^n\right][\zeta(n)-\eta(n+1)]
\mbox{ ~for ~} n\geq 3,\label{an}
\end{cases}\\
&{}&
\begin{cases}
b_0=5\log 2-\pi=0.324143\cdots, b_1=-\log
2+[35\zeta(3)-\pi^3]/8=0.690067\cdots,\\
b_n=5\lambda(2n+1)-\lambda(2n-1)-4\beta(2n+1) \mbox{ ~for ~} n\geq 2,\label{bn}
\end{cases}\\
&{}&~~c_n=1-\sum_{k=0}^nb_k \mbox{ ~and ~}
A_n=2^{n+1}\left(b_0-\sum_{k=0}^n2^{-k}a_k\right) \mbox{ ~for ~} n\in\IN\cup\{0\}.\label{cn}
\end{eqnarray}
Clearly, for $n\in\IN\cup\{0\}$,
\begin{eqnarray}
&{}&
\begin{cases}
a_{2n+2}=2[\zeta(2n+3)-\zeta(2n+1)]=-2\sum_{k=2}^{\infty}k^{-2n-3}(k^2-1)<0,\\
a_{2n+1}=2\left[\zeta(2n+1)-\left(1-2^{-2n-1}\right)\zeta(2n+2)\right]>4^{-n}\zeta(2n+2)>0,\label{SignOfan}
\end{cases}\\
&{}&~~~~c_{n+1}=c_n-b_{n+1}, ~A_{n+1}=2(A_n-a_{n+1}), ~n\in\IN\cup\{0\}.\label{(cn+1)+(An+1)}
\end{eqnarray}
Some properties of the constants $a_n$, $b_n$, $c_n$ and $A_n$ will be given in Lemma \ref{Lem1} and Corollaries \ref{Col1} and \ref{Col2}.
 We now state the main results of this paper.

\begin{theorem}\label{th1}
(1) The function $f(x)\equiv[1+x(1-x)]R(x)-B(x)$ has the following power series expansions
\begin{eqnarray}\label{SerOff(x)}
f(x)=\sum_{n=0}^{\infty}a_nx^n=\sum_{n=0}^{\infty}b_n(1-2x)^{2n},
~0<x\leq\frac12.
\end{eqnarray}

(2) $f$ ($f'$) is strictly decreasing and convex (increasing and convex) from $(0,1/2]$ onto
$[b_0,1)$ ($(a_1,0]$, respectively), and the function $F(x)\equiv b_0+b_1(1-2x)^2+B(x)-[1+x(1-x)]R(x)$ is strictly completely monotonic on $(0,1/2]$ with $F(0^+)=b_0+b_1-1=0.0142104\cdots$ and $F(1/2)=0$. Moreover, for each $n\in\IN\setminus\{1\}$, the function $F_{1,\thinspace n}(x)\equiv (-1)^{n+1}f^{(n)}(x)$ is strictly completely monotonic on $(0,1/2]$ with $F_{1,\thinspace n}(0^+)=(-1)^{n+1}n!a_n$, $F_{1,\thinspace n}(1/2)=0$ if $n$ is odd, and $F_{1,\thinspace n}(1/2)=-2^nn!b_{n/2}$ if $n$ is even.
In particular, for $x\in(0,1/2]$,
\begin{eqnarray}\label{Ineq1OfR}
\frac{b_0+(1-2x)P(x)+B(x)}{1+x(1-x)}\leq R(x)\leq\frac{b_0+(1-2x)Q(x)+B(x)}{1+x(1-x)},
\end{eqnarray}
with equality in each instance if and only if $x=1/2$, where $P(x)=\max\{0,1-b_0-b_1+b_1(1-2x)\}$ and $Q(x)=\min\{1-b_0,b_1(1-2x)\}$.

\end{theorem}

\begin{theorem}\label{th2}
For each $n\in\IN\cup\{0\}$ and $x\in(0,1/2]$, let $R_n(x)=\sum_{k=0}^na_kx^k$ and $S_n(x)=\sum_{k=0}^nb_k(1-2x)^{2k}$, $f$ as in Theorem 1.1, and define the functions $f_n$ and $g_n$ on $(0,1/2)$ by
\begin{eqnarray*}
f_n(x)=\frac{f(x)-R_n(x)}{x^{n+1}} \mbox{ ~and ~} g_n(x)=\frac{f(x)-S_n(x)}{(1-2x)^{2(n+1)}},
\end{eqnarray*}
respectively. Then we have the following conclusions:

(1) $f_0$ is strictly increasing and convex from $(0,1/2]$ onto
$(a_1,2(b_0-1)]$. Furthermore, for each $n\in\IN$, $f_{2n-1}$ ($f_{2n}$) is strictly increasing and concave (decreasing and
convex) from $(0,1/2]$ onto $(a_{2n},A_{2n-1}]$ ($[A_{2n},a_{2n+1})$, respectively). In particular, for each
$n\in\IN$ and all $x\in(0,1/2]$,
\begin{eqnarray}
\frac{B(x)+R_{2n+2}(x)+A_{2n+2}x^{2n+3}}{1+x(1-x)}&\leq& R(x)\leq
\frac{B(x)+R_{2n+1}(x)+A_{2n+1}x^{2(n+1)}}{1+x(1-x)},\label{Ineq2OfR}\\
\frac{B(x)+R_4(x)+A_4x^3}{1+x(1-x)}&\leq& R(x)\leq
\frac{B(x)+R_3(x)+A_3x^4}{1+x(1-x)},\label{Ineq3OfR}
\end{eqnarray}
with equality in each instance if and only if $x=1/2$.

(2) For each $n\in\IN\cup\{0\}$, $g_n$ is is strictly increasing and concave from $(0,1/2)$ onto $(c_n, b_{n+1})$. Furthermore, for each $m, n\in\IN\cup\{0\}$, the function $G_{n,\,m}(x)\equiv (-1)^{m+1}g_n^{(m)}(x)$ is strictly completely monotonic on $(0,1/2)$.
In particular, for each $n\in\IN\cup\{0\}$ and all $x\in(0,1/2]$,
\begin{eqnarray}
\frac{B(x)+S_{n+1}(x)+c_{n+1}(1-2x)^{2n+3}}{1+x(1-x)}&\leq&
R(x)\leq \frac{B(x)+S_{n+1}(x)}{1+x(1-x)},\label{Ineq4OfR}\\
\frac{B(x)+S_2(x)+c_2(1-2x)^5}{1+x(1-x)}&\leq& R(x)\leq
\frac{B(x)+S_2(x)}{1+x(1-x)},\label{Ineq5OfR}
\end{eqnarray}
with equality in each instance if and only if $x=1/2$.
\end{theorem}

By Theorem \ref{th1}, it is natural to ask whether the functions
\begin{eqnarray}\label{F1-F3}
F_1(x)\equiv
[1+x(1-x)]\frac{R(x)}{B(x)}, F_2(x)\equiv
\frac{R(x)}{B(x)}-\frac{1}{1+x(1-x)} \mbox{ ~and ~} F_3(x)\equiv R(x)-\frac{B(x)}{1+x(1-x)}
\end{eqnarray}
are monotone on $(0,1/2]$. Our next
result gives the answer to this question.

\begin{theorem}\label{th3}
(1) There exists a unique number
$x_1\in(1/4,1/2)$ such that the function $F_1$ defined by (\ref{F1-F3}) is strictly increasing on $(0,x_1]$, and decreasing on $[x_1,1/2]$, with $F_1(0^+)=1$ and $F_1(1/2)=(5\log2)/\pi=1.103178\cdots$. However, $F_1$ is neither convex nor concave on $(0,1/2]$. In particular, for $x\in(0,1/2]$,
\begin{eqnarray}\label{Ineq6OfR}
\frac{\alpha B(x)}{1+x(1-x)}<R(x)\leq \frac{\delta
B(x)}{1+x(1-x)},
\end{eqnarray}
with the best possible coefficients $\alpha=1$ and $\delta=F_1(x_1)$. Moreover,
\begin{eqnarray}\label{delta}
1.111592\cdots=\frac{19\sqrt{2}\log8}{16\pi}<\delta<
1.112146.
\end{eqnarray}

(2) There exists a unique number $x_2\in(0,1/2)$ such that the
function $F_2$ defined by (\ref{F1-F3}) is strictly
increasing on $(0,x_2]$, and decreasing on $[x_2,1/2]$ with
$F_2(0^+)=0$ and $F_2(1/2)=4b_0/(5\pi)=0.082542\cdots$. However, $F_2$ is neither convex nor concave on $(0,1/2]$.

(3) The function $F_3$ defined by (\ref{F1-F3}) is
strictly decreasing and convex from $(0,1/2]$ onto $[\rho,1)$,
where $\rho=4b_0/5=\log16-4\pi/5=0.259314\cdots$. In particular, for $x\in(0,1/2]$,
\begin{eqnarray}\label{Ineq7OfR}
\rho+\frac{B(x)}{1+x(1-x)}\leq R(x)\leq
\rho+(1-\rho)(1-2x)+\frac{B(x)}{1+x(1-x)},
\end{eqnarray}
with equality in each instance if and only if $x=1/2$.
\end{theorem}

\section{\normalsize Preliminaries}

In this section, we prove two technical lemmas needed in the proofs of our main results
stated in Section 1. Our first lemma shows some properties of $b_n$, $\lambda(n)$ and $\beta(n)$.

\begin{lemma}\label{Lem1}
(1) The functions
$\lambda(x)\equiv\sum_{k=0}^{\infty}(2k+1)^{-x}$ and $\varphi_1(x)\equiv \lambda(x)-\lambda(x+1)$ are both strictly
decreasing and convex on $(1,\infty)$, with $\lambda((1,\infty))=(1,\infty)$ and $\varphi_1((1,\infty))=(0,\infty)$. Moreover, for each $c\in[C_1,\infty)$, the function $\varphi_2(x)\equiv \lambda(x+c)/\lambda(x)$
is strictly increasing from $[2,\infty)$ onto $[C_2,1)$, where $C_1=\log(\pi^2/8)]/\log3=0.191166\cdots$ and $C_2=8\lambda(2+c)/\pi^2>8/\pi^2=0.810569\cdots$.

(2) The function $\beta(x)\equiv
\sum_{k=0}^{\infty}(-1)^k(2k+1)^{-x}$ is strictly increasing from
$[1,\infty)$ onto $[\pi/4,1)$, and concave on $[2/\log3,\infty)$.

(3) Let $\varphi_3(x)=5\lambda(x+1)-\lambda(x-1)-4\beta(x+1)$, $\varphi_4(x)=(x-3+1/\log5)\varphi_3(x)$, $\varphi_5(x)=x\varphi_3(x)$ and $\varphi_6(x)=\varphi_5(x+2)-\varphi_5(x)$, for $x\in[3,\infty)$. Then $\varphi_3$ is strictly increasing and
concave from $[3,\infty)$ onto $[\mu_1,0)$, $\varphi_4$ and $\varphi_5$ are both
strictly increasing on $[3,\infty)$ with $\varphi_4([3,\infty))=[\mu_2,0)$ and $\varphi_5([3,\infty))=[\mu_3,0)$, and $\varphi_6$ is strictly decreasing from $[3,\infty)$ onto $(0, \mu_4]$, where $\mu_1=\varphi_3(3)=(5\pi^4/96)-4\beta(4)-\pi^2/8=-0.116088\cdots$, $\mu_2=\mu_1/\log5=-0.072129\cdots$, $\mu_3=3\mu_1=-0.348265\cdots$ and $\mu_4=\varphi_5(5)-\mu_3=0.337348\cdots$.

(4) For $n\in\IN$, the sequences $\{b_{n+1}\}$ and $\{(n+1)b_{n+1}\}$ are both strictly increasing, while the sequence $\{(n+2)b_{n+2}-(n+1)b_{n+1}\}$ is strictly decreasing, with $\lim_{n\to\infty}nb_n=\lim_{n\to\infty}[(n+1)b_{n+1}-nb_n ]=0$. In particular, for $n\in\IN\setminus\{1\}$,
\begin{eqnarray}\label{IneqOfbn}
-0.027624\cdots=2b_2\leq nb_n<(n+1)b_{n+1}\leq \min\{0,nb_n+3b_3-2b_2\},
\end{eqnarray}
with equality in each instance if and only if $n=2$.
\end{lemma}

{\it Proof.} (1) The monotonicity of $\lambda$ is clear. Since the derivative
\begin{eqnarray}\label{lam'}
\lambda'(x)=-\sum_{k=1}^{\infty}\frac{\log(2k+1)}{(2k+1)^x}
\end{eqnarray}
is clearly increasing on $(1,\infty)$, the function
$\lambda$ is convex on $(1,\infty)$.

It is easy to see that $\varphi_1(1)=\infty$ and $\varphi_1(\infty)=0$. By (\ref{lam'}), we have
\begin{eqnarray}\label{lam'(x+1)-lam'}
\varphi_1'(x)=-2\sum_{k=1}^{\infty}\frac{k}{(2k+1)^{x+1}}\log(2k+1),
\end{eqnarray}
which is negative and strictly increasing on $(1,\infty)$, and hence the result for $\varphi_1$ follows.

Clearly, $\varphi_2(2)=\lambda(2+c)/\lambda(2)=8\lambda(2+c)/\pi^2=C_2$ and
$\varphi_2(\infty)=1$. By differentiation and (\ref{lam'}),
\begin{eqnarray}\label{DerOfphi2}
\lambda(x)^2\varphi_2'(x)
=\sum_{k=1}^{\infty}\frac{\log(2k+1)}{(2k+1)^{x+c}}\left[(2k+1)^c\lambda(x+c)-\lambda(x)\right]
>\varphi_7(x)\sum_{k=1}^{\infty}\frac{\log(2k+1)}{(2k+1)^{x+c}},
\end{eqnarray}
where $\varphi_7(x)=3^c\lambda(x+c)-\lambda(x)$, and by (\ref{lam-beta}),
\begin{eqnarray*}
\varphi_7(x)=\sum_{k=0}^{\infty}\frac{3^c-(2k+1)^c}{(2k+1)^{x+c}}
=3^c-1-\sum_{k=1}^{\infty}\frac{(2k+1)^c-3^c}{(2k+1)^{x+c}}
\end{eqnarray*}
which is clearly strictly increasing on $(1,\infty)$. Since $c\geq [\log(\pi^2/8)]/\log3=[\log\lambda(2)]/\log3$,
$$\varphi_7(2)=3^c\lambda(2+c)-\lambda(2)\geq \lambda(2)[\lambda(2+c)-1]>0.$$
Hence the result for $\varphi_2$ follows from (\ref{DerOfphi2}).

(2) Clearly,
$\lim_{x\to \infty}\beta(x)=1$ and $\beta(1)=\pi/4$. Set $\beta_1(t)=(\log t)/t$ for $t\in[3,\infty)$. Then by differentiation,
\begin{eqnarray}\label{DerOfbeta}
\beta^{\thinspace\prime}(x)&=&\sum_{k=1}^{\infty}(-1)^{k+1}\frac{\log(2k+1)}{(2k+1)^x}
=\sum_{k=1}^{\infty}\frac{\log[2(2k-1)+1]}{[2(2k-1)+1]^x}-\sum_{k=1}^{\infty}\frac{\log[2(2k)+1]}{[2(2k)+1]^x}\nonumber\\
&=&\frac{1}{x}\sum_{k=1}^{\infty}\left[\beta_1((4k-1)^x)-\beta_1((4k+1)^x)\right].
\end{eqnarray}
It is easy to verify that the function $\beta_1$ is strictly decreasing on
$[\e,\infty)$. Since $(4k-1)^x\geq 3^x\geq 3>\e$ for $k\geq1$ and
$x\geq1$, $\beta^{\thinspace\prime}(x)>0$ for $x\in[1,\infty)$ by (\ref{DerOfbeta}), so
that the monotonicity of $\beta(x)$ follows.

By the first equality in (\ref{DerOfbeta}) and by differentiation,
\begin{eqnarray*}
\beta^{\thinspace\prime\prime}(x)=\frac{1}{x^2}\sum_{k=1}^{\infty}\left[\beta_2((4k+1)^x)-\beta_2((4k-1)^x)\right],
\end{eqnarray*}
where $\beta_2(t)=(\log t)^2/t$ for $t\geq 3$. It is easy to show
that $\beta_2$ is strictly decreasing on $[\e^2,\infty)$. Hence if
$(4k-1)^x\geq 3^x\geq \e^2$, that is, $x\geq2/\log3$, then
$\beta^{\thinspace\prime\prime}(x)<0$. This shows that $\beta$ is concave on
$[2/\log3,\infty)$.

(3) It follows from (\ref{zeta(2)}) and [1,Table
23.3] that
$$\varphi_3(3)=\mu_1=5\lambda(4)-\lambda(2)-4\beta(4)=\frac{5\pi^4}{96}-\frac{\pi^2}{8}-4\beta(4)=-0.116088\cdots.$$
Clearly, $\varphi_3(\infty)=0$, and $k(k+1)-1+(-1)^k\geq 6$ for all $k\geq2$. By (\ref{lam-beta}), we have
\begin{eqnarray}
\varphi_3(x)&=&-4\sum_{k=2}^{\infty}\frac{k(k+1)-1+(-1)^k}{(2k+1)^{x+1}},\label{phi3}\\
\varphi_3'(x)&=&4\sum_{k=2}^{\infty}\frac{k(k+1)-1+(-1)^k}{(2k+1)^{x+1}}\log(2k+1).\label{phi3'}
\end{eqnarray}
Hence the result for $\varphi_3$ follow from (\ref{phi3}) and (\ref{phi3'}).

Next, for each $k\in\IN\setminus\{1\}$ and for $x\in[3,\infty)$, let
$\varphi_{8,\,k}(x)=(x-3+1/\log5)(2k+1)^{-x}$. Then $\varphi_{8,\,k}(3)=1/[(2k+1)^3\log5]$, $\varphi_{8,\,k}(\infty)=0$,
and by (\ref{phi3}),
\begin{eqnarray}\label{phi4}
\varphi_4(x)=-4\sum_{k=2}^{\infty}\frac{k(k+1)-1+(-1)^k}{2k+1}\varphi_{8,\,k}(x)
\end{eqnarray}
Clearly, $\varphi_4(3)=\mu_2$, and $\varphi_4(\infty)=0$ by (\ref{phi4}). Since
\begin{eqnarray*}
\varphi_{8,\,k}'(x)&=&(2k+1)^{-x}[1-(x-3+1/\log5)\log(2k+1)]\\
&\leq& -(x-3)(2k+1)^{-x}\log5<0 \mbox{ ~for ~} x>3,
\end{eqnarray*}
$\varphi_{8,\,k}$ is strictly decreasing from
$[3,\infty)$ onto $(0,\varphi_{8,\,k}(3)]$, so that the monotonicity of
$\varphi_4$ follows from (\ref{phi4}).

Since $3-1/\log5=2.378665\cdots>0$, and since
$$-\varphi_5(x)=[-\varphi_4(x)][1+(3-1/\log5)/(x-3+1/\log5)]$$
which is a product of two positive and strictly decreasing functions on $[3,\infty)$, the monotonicity of $\varphi_5$ follows.

Clearly, $\varphi_6(3)= \mu_4=\varphi_5(5)-\mu_3=5[5\lambda(6)-\lambda(4)-4\beta(6)]-\mu_3=0.337348\cdots$ and $\varphi_6(\infty)=0$. It follows from (\ref{phi3}) that
\begin{eqnarray}\label{phi6}
\varphi_6(x)=(x+2)\varphi_3(x+2)-x\varphi_3(x)=8\sum_{k=2}^{\infty}\frac{k(k+1)-1+(-1)^k}{(2k+1)^{x+3}}[2k(k+1)x-1].
\end{eqnarray}
Differentiation gives
\begin{eqnarray*}
\frac{d}{dx}\left[\frac{2\left(k^2+k\right)x-1}{(2k+1)^{x+3}}\right]
&=&\frac{2k(k+1)\log(2k+1)}{(2k+1)^{x+3}}\left[\frac{1}{\log(2k+1)}+\frac{1}{2k(k+1)}-x\right]\\
&\leq&\frac{2k(k+1)\log(2k+1)}{(2k+1)^{x+3}}\left(\frac{1}{\log5}+\frac{1}{12}-3\right)<0
\end{eqnarray*}
for $x\in[3,\infty)$, since $(1/\log5)+(1/12)=0.704668\cdots$. Hence the monotonicity of $\varphi_6$ follows from (\ref{phi6}).

(4) Clearly, $b_{n+1}=\varphi_3(2n+2)$, $(n+1)b_{n+1}=\varphi_5(2n+2)/2$ and $(n+2)b_{n+2}-(n+1)b_{n+1}=\varphi_6(2n+2)/2$. Hence part (4) follows from part (3). $\Box$

\begin{lemma}\label{Lem2}
For $x\in(0,1/2]$,
\begin{eqnarray}
R(x)&=&\frac{1}{x}+\sum_{n=1}^{\infty}\left[1+(-1)^n\right]\zeta(n+1)x^n =\log
16+4\sum_{n=1}^{\infty}\lambda(2n+1)(1-2x)^{2n},\label{SerOfR}\\
B(x)&=&\frac{1}{x}+\sum_{n=1}^{\infty}\left[1+(-1)^{n+1}\right]\eta(n+1)x^n
=4\sum_{n=0}^{\infty}\beta(2n+1)(1-2x)^{2n},\label{SerOfB}\\
R\left(\frac14\right)&=&-2\gamma-\psi\left(\frac14\right)-\psi\left(\frac34\right)=6\log2=4.158883\cdots,\label{R(1/4)}\\
R'\left(\frac14\right)&=&\psi'\left(\frac34\right)-\psi'\left(\frac14\right)=-16\beta(2)=-14.655449\cdots.\label{R'(1/4)}
\end{eqnarray}

\end{lemma}

{\it Proof.} It is well known that
\begin{eqnarray}
\psi(1+x)&=&\psi(x)+\frac{1}{x},\label{psi1}\\
\psi(1+x)+\gamma&=&\sum_{n=1}^{\infty}(-1)^{n+1}\zeta(n+1)x^n, ~|x|<1,\label{psi2}\\
\psi^{(n)}(x)&=&(-1)^{n+1}n!\sum_{k=0}^{\infty}\frac{1}{(k+x)^{n+1}}.\label{n-Dpsi}
\end{eqnarray}
(Cf. \citep[6.3.5, 6.3.14 \& 6.4.10]{1}.) By (\ref{psi1}) and (\ref{psi2}),
\begin{eqnarray*}
R(x)&=&\frac{1}{x}-[\gamma+\psi(1-x)]-[\gamma+\psi(1+x)]\\
&=&\frac{1}{x}+\sum_{n=1}^{\infty}(-1)^n\zeta(n+1)(-x)^n+\sum_{n=1}^{\infty}(-1)^n\zeta(n+1)x^n,
\end{eqnarray*}
yielding the first equality in (\ref{SerOfR}). By differentiation,
\begin{eqnarray}\label{n-DR}
R^{(n)}(x)=(-1)^{n+1}\psi^{(n)}(1-x)-\psi^{(n)}(x).
\end{eqnarray}
 Hence by (\ref{n-Dpsi}) and (\ref{n-DR}), $R^{(2n-1)}(1/2)=0$ and
\begin{eqnarray*}
R^{(2n)}\left(\frac12\right)=-2\psi^{(2n)}\left(\frac12\right)=\sum_{k=0}^{\infty}\frac{2(2n)!}{(k+1/2)^{2n+1}}=4^{n+1}(2n)!\lambda(2n+1),
\end{eqnarray*}
so that $R(x)$ has the following power series
expansion
\begin{eqnarray*}
R(x)&=&R\left(\frac12\right)+\sum_{n=1}^{\infty}\frac{R^{(n)}(1/2)}{n!}\left(x-\frac{1}{2}\right)^n
=\log16+\sum_{n=1}^{\infty}\frac{R^{(2n)}(1/2)}{(2n)!}\left(x-\frac{1}{2}\right)^{2n}\\
&=&\log16+4\sum_{n=1}^{\infty}\lambda(2n+1)(1-2x)^{2n}
\end{eqnarray*}
at $x=1/2$. This yields the second equality in (\ref{SerOfR}).

By \citep[4.3.68 \& 23.1.18]{1},
\begin{eqnarray}\label{SerOfSin}
\frac{1}{\sin
z}=\frac{1}{z}+2\sum_{n=1}^{\infty}\frac{\eta(2n)}{\pi^{2n}}z^{2n-1}
=\frac{1}{z}+\sum_{n=1}^{\infty}\left[1+(-1)^{n+1}\right]\frac{\eta(n+1)}{\pi^{n+1}}z^{n},
\end{eqnarray}
for $|z|<\pi$. Hence the first equality in (\ref{SerOfB}) follows from (\ref{B})
and (\ref{SerOfSin}). By \citep[4.3.69 \& 23.1.18]{1}, we have
\begin{eqnarray}\label{SerOfCos}
\frac{1}{\cos
z}=2\sum_{n=0}^{\infty}\left(\frac{2}{\pi}\right)^{2n+1}\beta(2n+1)z^{2n}, ~|z|<\frac{\pi}{2},
\end{eqnarray}
from which it follows that
$$B(x)=\frac{\pi}{\sin(\pi x)}=\frac{\pi}{\cos[\pi(1-2x)/2]}
=4\sum_{n=0}^{\infty}\beta(2n+1)(1-2x)^{2n},$$
yielding the
second equality in (\ref{SerOfB}).

Next, by \citep[6.3.3 \& 6.3.8]{1},
\begin{eqnarray}\label{psi(2x)}
\psi\left(\frac12\right)=-\gamma-\log4, ~\psi(2x)=\frac12\psi(x)+\frac12\psi\left(x+\frac12\right)+\log 2,
\end{eqnarray}
from which we obtain
\begin{eqnarray*}
R\left(\frac14\right)=-2\gamma-\left[\psi\left(\frac14\right)+\psi\left(\frac34\right)\right]
=-2\gamma-2\left[\psi\left(\frac12\right)-\log2\right]=6\log2.
\end{eqnarray*}

Finally, by (\ref{lam-beta}), it is easy to see that
\begin{eqnarray}\label{beta(2)}
\beta(2)=\sum_{n=0}^{\infty}\frac{(-1)^n}{(2n+1)^2}=\sum_{n=0}^{\infty}\frac{1}{[2(2n)+1]^2}-\sum_{n=0}^{\infty}\frac{1}{[2(2n+1)+1]^2}.
\end{eqnarray}
It follows from (\ref{n-Dpsi}), (\ref{beta(2)}) and \citep[Table 23.3]{1} that
\begin{eqnarray*}
R'\left(\frac14\right)=\psi'\left(\frac34\right)-\psi'\left(\frac14\right)
=16\left[\sum_{n=0}^{\infty}\frac{1}{(4n+3)^2}-\sum_{n=0}^{\infty}\frac{1}{(4n+1)^2}\right]
=-16\beta(2). ~\Box
\end{eqnarray*}

\begin{remark}\label{Rem1}
In \citep[Theorem 2.2]{16}, it was proved that
$R(x)=(1/x)+2\sum_{n=1}^{\infty}\zeta(2n+1)x^{2n}$ for $x\in(0,1/2]$,
which is consistent with the first equality in (\ref{SerOfR}). However, its proof given in \cite{16} is quite
complicated.
\end{remark}

\section{\normalsize Proof of Theorem \ref{th1}}

(1) Since $f(x)=R(x)-B(x)+x(1-x)R(x)$, it follows from Lemma \ref{Lem2},
(\ref{lambda(n+1)}) and (\ref{an}) that
\begin{eqnarray*}
f(x)&=&\sum_{n=1}^{\infty}\left\{\left[1+(-1)^n\right]\zeta(n+1)-\left[1+(-1)^{n+1}\right]\eta(n+1)\right\}x^n\\
&{}&+(1-x)\left\{1+\sum_{n=1}^{\infty}\left[1+(-1)^n\right]\zeta(n+1)x^{n+1}\right\}\\
&=&1-x+\sum_{n=1}^{\infty}\left\{\left[1+(-1)^n\right]\zeta(n+1)-\left[1+(-1)^{n+1}\right]\eta(n+1)\right\}x^n\\
&{}&+\sum_{n=1}^{\infty}[1+(-1)^n]\zeta(n+1)x^{n+1}-\sum_{n=1}^{\infty}[1+(-1)^n]\zeta(n+1)x^{n+2}\\
&=&1-[1+2\eta(2)]x+2\zeta(3)x^2+\sum_{n=3}^{\infty}\left\{\left[1+(-1)^n\right]\zeta(n+1)-\left[1+(-1)^{n+1}\right]\eta(n+1)\right\}x^n\\
&{}&+\sum_{n=2}^{\infty}\left[1+(-1)^{n+1}\right]\zeta(n)x^n-\sum_{n=3}^{\infty}\left[1+(-1)^n\right]\zeta(n-1)x^n\\
&=&1-[1+2\eta(2)]x+2\zeta(3)x^2+\sum_{n=3}^{\infty}\left\{\left[1+(-1)^n\right]\zeta(n+1)+\left[1+(-1)^{n+1}\right]\zeta(n)\right.\\
&{}&\left.-\left[1+(-1)^{n+1}\right]\eta(n+1)-\left[1+(-1)^n\right]\zeta(n-1)\right\}x^n\\
&=&1-[1+\zeta(2)]x+2\zeta(3)x^2+\sum_{n=3}^{\infty}a_nx^n=\sum_{n=0}^{\infty}a_nx^n,
\end{eqnarray*}
which yields the first power series expansion in (\ref{SerOff(x)}).

Clearly, $1+x(1-x)=[5-(1-2x)^2]/4$. It follows from (\ref{SerOfR}),
(\ref{SerOfB}), (\ref{lambda(n+1)}), (\ref{zeta(2)}) and (\ref{bn}) that
\begin{eqnarray*}
f(x)&=&\left[\frac{5}{4}-\frac{(1-2x)^2}{4}\right]\left[\log 16
+4\sum_{n=1}^{\infty}\lambda(2n+1)(1-2x)^{2n}\right]-4\sum_{n=0}^{\infty}\beta(2n+1)(1-2x)^{2n}\\
&=&b_0-(1-2x)^2\log2-\sum_{n=1}^{\infty}\lambda(2n+1)(1-2x)^{2n+2}+\sum_{n=1}^{\infty}[5\lambda(2n+1)-4\beta(2n+1)](1-2x)^{2n}\\
&=&b_0+[5\lambda(3)-4\beta(3)-\log2](1-2x)^2+\sum_{n=2}^{\infty}[5\lambda(2n+1)-\lambda(2n-1)-4\beta(2n+1)](1-2x)^{2n}\\
&=&b_0+b_1(1-2x)^2+\sum_{n=2}^{\infty}b_n(1-2x)^{2n}=\sum_{n=0}^{\infty}b_n(1-2x)^{2n},
\end{eqnarray*}
yielding the second series expansion in (\ref{SerOff(x)}).

(2) Clearly, $f(0^+)=a_0=1$ and $f(1/2)=b_0$ by (\ref{SerOff(x)}). By (\ref{SerOff(x)}) and differentiation, we
obtain
\begin{eqnarray}
f'(x)&=&-2\sum_{k=1}^{\infty}(2k)b_k(1-2x)^{2k-1},\label{SerOff'}\\
f''(x)&=&(-2)^2\sum_{k=0}^{\infty}(2k+2)(2k+1)b_{k+1}(1-2x)^{2k},\label{SerOff''}\\
f'''(x)&=&(-2)^3\sum_{k=1}^{\infty}(2k+2)(2k+1)(2k)b_{k+1}(1-2x)^{2k-1}.\label{SerOff'''}
\end{eqnarray}
Generally, by the mathematical induction, it is not difficult to
obtain the following two expressions
\begin{eqnarray}
f^{(2n-1)}(x)&=&-2^{2n-1}\sum_{k=1}^{\infty}\frac{(2n+2k-2)!}{(2k-1)!})b_{k+n-1}(1-2x)^{2k-1},\label{SerOf(2n-1)-Df}\\
f^{(2n)}(x)&=&2^{2n}\sum_{k=0}^{\infty}\frac{(2n+2k)!}{(2k)!}b_{n+k}(1-2x)^{2k},\label{SerOf(2n)-Df}
\end{eqnarray}
for $n\in\IN$. By (\ref{SerOf(2n-1)-Df}) and (\ref{SerOf(2n)-Df}), $f^{(2n-1)}(1/2)=0$ and
$f^{(2n)}(1/2)=4^n(2n)!b_n$ for $n\in\IN$. The limiting value
$f^{(n)}(0^+)=n!a_n$ follows from the first equality in (\ref{SerOff(x)}). Hence it follows from (\ref{IneqOfbn}),
(\ref{SerOf(2n-1)-Df}) and (\ref{SerOf(2n)-Df}) that for all $n\geq 2$, $f^{(2n-1)}$ (
$f^{(2n)}$ ) is strictly decreasing and convex (increasing and
concave) from $(0,1/2]$ onto $[0,(2n-1)!a_{2n-1})$ (
$((2n)!a_{2n},4^n(2n)!b_n]$, respectively).

It follows from (\ref{IneqOfbn}) and (\ref{SerOff''}) that $f''$ is strictly increasing and concave on $(0,1/2]$,
with $f''(0^+)=2a_2>0$ by (\ref{SerOff(x)}). Hence $f'$ is strictly increasing and
convex on $(0,1/2]$ with $f'(1/2)=0$ by (\ref{SerOff'}) and $f'(0^+)=a_1$ by (\ref{SerOff(x)}), and $f$ is strictly
decreasing and convex from $(0,1/2]$ onto $[b_0,1)$.

Clearly, $F(0^+)=b_0+b_1-f(0^+)=b_0+b_1-1=0.0142104\cdots$ and $F(1/2)=b_0-f(1/2)=0$.
 It follows from (\ref{SerOff(x)}) that $F$ has the following power series expansion
\begin{eqnarray}\label{SerOfF(x)}
F(x)=-\sum_{k=2}^{\infty}b_k(1-2x)^{2k}.
\end{eqnarray}
Applying Lemma \ref{Lem1}(4) and (\ref{SerOfF(x)}), one can easily see that
$F$ is strictly completely monotonic on $(0,1/2]$.

Next, for $m, n, k\in\IN$ and $x\in(0,1/2)$,
\begin{eqnarray}\label{F1n}
(-1)^mF_{1,\,n}^{(m)}(x)=(-1)^{m+n+1}f^{(n+m)}(x)=
\begin{cases}
-f^{(2k)}(x), \mbox{~~if ~} m+n=2k,\\
f^{(2k-1)}(x), \mbox{~if ~} m+n=2k-1.\\
\end{cases}
\end{eqnarray}
If $n\geq2$, then $m+n\geq3$ and $k\geq2$. It follows from Lemma \ref{Lem1}(4) and the monotonicity properties of $f^{(2k-1)}$ and
$f^{(2k)}$, $k\geq 2$, that for $x\in(0,1/2)$,
\begin{eqnarray*}
-f^{(2k)}(x)&>&-f^{(2k)}(1/2)=-4^k(2k)!b_{k}>0,\\
f^{(2k-1)}(x)&>&f^{(2k-1)}(1/2)=0.
\end{eqnarray*}
Hence by (\ref{F1n}), $(-1)^mF_{1,\,n}^{(m)}(x)>0$ for each $m\in\IN$, $n\in\IN\setminus\{1\}$ and for all $x\in(0,1/2)$. This yields the strictly complete monotonicity of $F_{1,\,n}$.

Clearly, $F_{1,\thinspace n}(0^+)=(-1)^{n+1}f^{(n)}(0^+)=(-1)^{n+1}n!a_n$. By (\ref{SerOf(2n-1)-Df}) and (\ref{SerOf(2n)-Df}), $F_{1,\thinspace n}(1/2)=0$ if $n$ is odd, and $F_{1,\thinspace n}(1/2)=-2^nn!b_{n/2}$ if $n$ is even.

Finally, it follows from the monotonicity and convexity properties of $f$ and $F$ that
\begin{eqnarray}
\frac{b_0+B(x)+(1-b_0)(1-2x)}{1+x(1-x)}&\geq& R(x)\geq \frac{b_0+B(x)}{1+x(1-x)}
,\label{Ineq0OfR}\\
\frac{b_0+b_1(1-2x)^2+B(x)}{1+x(1-x)}&\geq& R(x)\geq \frac{b_0+(1-b_0-b_1)(1-2x)+b_1(1-2x)^2+B(x)}{1+x(1-x)}
,\label{Ineq00OfR}
\end{eqnarray}
with equality in each instance if and only if $x=1/2$. Hence (\ref{Ineq1OfR}) and its equality case follow. $\Box$

\begin{corollary}\label{Col1}
Let $a_n$,
$b_n$ and $c_n$ be as in (\ref{an})--(\ref{cn}), and set $d_n=\sum_{k=0}^nb_k\lambda(2n-2k+2)-(n+1)b_{n+1}$ for $n\in\IN\cup\{0\}$, $D_0=5d_0-b_0$ and $D_n=5d_n-d_{n-1}-b_n$ for $n\in\IN$. Then we have the following conclusions:

(1) The following identities hold
\begin{eqnarray}
&{}&\sum_{k=0}^{\infty}a_k=\sum_{k=0}^{\infty}b_k=1,\label{Id1}\\
&{}&\sum_{k=0}^{\infty}2^{-k}a_k=b_0=5\log2-\pi.\label{Id2}
\end{eqnarray}

(2) For $n\in\IN$, the sequence $\{c_n\}$ is strictly increasing. In particular, for $n\in\IN\setminus\{1\}$,
\begin{eqnarray}
-0.014210\cdots=c_1&<&-0.000398\cdots=c_2\leq c_n<\lim_{n\to\infty} c_n=0,\label{Ineq1Ofcn}\\
1&<&\sum_{k=0}^nb_k=1-c_n\leq1-c_2.\label{IneqOf(1-cn)}
\end{eqnarray}

(3) $d_0=-0.290171\cdots$, $d_1=1.207861\cdots$, $d_2=1.008920\cdots$, $d_3=1.000824\cdots$, $\lim_{n\to\infty}d_n=1$, and
\begin{eqnarray}\label{IneqOfdn}
d<d_n<\widetilde{d}<d_2 \mbox{ ~for ~} n\geq3,
\end{eqnarray}
where $d=(1-\pi^2/8)(b_0+b_1)+\pi^2/8=0.996679\cdots$ and $\widetilde d=d_3-[(\pi^4/96)-1]b_2-[(\pi^2/8)-1]b_3=1.001117\cdots$.

(4) $D_0=[(5\pi^2/8)-1]b_0-5b_1=-1.7750006\cdots$, $D_1=5d_1-d_0-b_1=5.639413\cdots$, $D_2=5d_2-d_1-b_2=3.850551\cdots$, $D_3=5d_3-d_2-b_3=3.995587\cdots$, and
\begin{eqnarray}\label{Dn}
D_n>5d-\widetilde{d}=3.982277\cdots, ~n\geq4.
\end{eqnarray}
\end{corollary}

\smallskip
{\bf Proof.} (1) Let $f$ be as in Theorem \ref{th1}. Then $f(1-x)=f(x)$ by the symmetry. Hence $f(1^-)=f(0^+)=a_0=1$, so that $\sum_{k=0}^{\infty}a_k=\sum_{k=0}^{\infty}b_k=1$. We obtain (\ref{Id2}) by taking $x=1/2$ in (\ref{SerOff(x)}).

(2) By (\ref{Id1}) and (\ref{cn}), $\lim_{n\to\infty}c_n=0$. Since $c_{n+1}-c_n=-b_{n+1}>0$ for $n\in\IN$ by Lemma \ref{Lem1}(4) and (\ref{(cn+1)+(An+1)}), the monotonicity of $c_n$ and (\ref{Ineq1Ofcn}) follow. (\ref{IneqOf(1-cn)}) holds by (\ref{Ineq1Ofcn}).

(3) By (\ref{lambda(n+1)}), (\ref{zeta(2)}), (\ref{bn}) and \citep[Table 23.3]{1}, we obtain
\begin{eqnarray*}
d_0&=&b_0\lambda(2)-b_1=(\pi^2b_0/8)-b_1=-0.290171\cdots,\\
d_1&=&b_0\lambda(4)+b_1\lambda(2)-2[5\lambda(5)-\lambda(3)-4\beta(5)]=1.207861\cdots,\\
d_2&=&b_0\lambda(6)+b_1\lambda(4)+b_2\lambda(2)-3b_3=1.008920\cdots,\\
d_3&=&b_0\lambda(8)+b_1\lambda(6)+b_2\lambda(4)+b_3\lambda(2)-4b_4=1.000824\cdots.
\end{eqnarray*}

Applying Lemma \ref{Lem1}(1) and (4), one can easily show that $\lim_{n\to\infty}\sum_{k=2}^nb_k[\lambda(2n-2k+2)-1]=0$ by the definition of limit. Hence by Lemma \ref{Lem1}(1) and (4),
\begin{eqnarray*}
\lim_{n\to\infty}d_n&=&\lim_{n\to\infty}\left\{b_0\lambda(2n+2)+b_1\lambda(2n)
+\sum_{k=2}^nb_k-(n+1)b_{n+1}\right.\\
&{}&\left.+\sum_{k=2}^nb_k[\lambda(2n-2k+2)-1]\right\}=\sum_{k=0}^{\infty}b_k=1.
\end{eqnarray*}

Since $\lambda(m)>1$ and $b_m<0$ for $m\geq2$ by Lemma \ref{Lem1}(1) and (4), it follows that for $n\geq3$,
\begin{eqnarray*}
d_n&=&b_0\lambda(2n+2)+b_1\lambda(2n)+\sum_{k=2}^nb_k\lambda(2n-2k+2)-(n+1)b_{n+1}\\
&<&\widetilde{d_n}\equiv b_0\lambda(2n+2)+b_1\lambda(2n)+\sum_{k=2}^nb_k-(n+1)b_{n+1}\\
&=& b_0[\lambda(2n+2)-1]+b_1[\lambda(2n)-1]+1-c_n-(n+1)b_{n+1},
\end{eqnarray*}
which strictly decreasing in $n$ by part (2) and Lemma \ref{Lem1}(1) and (4). Hence for $n\geq3$,
\begin{eqnarray*}
\widetilde{d_n}&\leq&\widetilde{d_3}=\widetilde{d}\equiv b_0[\lambda(8)-1]+b_1[\lambda(6)-1]+1-c_3-4b_4\\
&=&d_3-[\lambda(4)-1]b_2-[\lambda(2)-1]b_3=1.001117\cdots<d_2.
\end{eqnarray*}
Hence the second and third inequalities in (\ref{IneqOfdn}) hold.

On the other hand, by part (2), and by Lemma \ref{Lem1}(1) and (4), we have
\begin{eqnarray*}
d_n&\geq&b_0\lambda(2n+2)+b_1\lambda(2n)+\lambda(2)\sum_{k=2}^nb_k-(n+1)b_{n+1}=\overline{d_n}\\
&\equiv& b_0[\lambda(2n+2)-\lambda(2)]+b_1[\lambda(2n)-\lambda(2)]+\lambda(2)(1-c_n)-(n+1)b_{n+1},
\end{eqnarray*}
which strictly decreasing in $n$ by part (2) and Lemma \ref{Lem1}(1) and (4), and hence
\begin{eqnarray*}
d_n\geq\overline{d_n}>\lim_{n\to\infty}\overline{d_n}=d\equiv b_0\left(1-\frac{\pi^2}{8}\right)+b_1\left(1-\frac{\pi^2}{8}\right)+\frac{\pi^2}{8}
=0.996679\cdots.
\end{eqnarray*}
This yields the first inequality in (\ref{IneqOfdn}).

(4) By computation, one can obtain the values of $D_n$ for $0\leq n\leq 3$. It follows from  (\ref{IneqOfbn}) and (\ref{IneqOfdn}) that
\begin{eqnarray*}
D_n=5d_n-d_{n-1}-b_n>5d-\widetilde{d}=3.982277\cdots
\end{eqnarray*}
for $n\geq 4$, which yields the inequality (\ref{Dn}). $\Box$

\section{\normalsize Proof of Theorem \ref{th2}}

(1) Clearly, $f_0(0^+)=f'(0^+)=a_1$ and $f_0(1/2)=2(b_0-1)$. Let $h_0(x)=xf'(x)-[f(x)-a_0]$. Then
\begin{eqnarray*}
h_0(0^+)=0, ~f_0'(x)=\frac{h_0(x)}{x^2}, ~h_0'(x)\left[\frac{d}{dx}(x^2)\right]^{-1}=\frac12 f''(x).
\end{eqnarray*}
By Theorem \ref{th1}(2), $f''$ is strictly increasing on $(0,1/2]$, and hence so is $f_0'$ by \citep[Theorem 1.25]{2}. Since $f_0'(0^+)=f''(0^+)/2=a_2>0$ by l'H\^opital's Rule, the monotonicity and convexity properties of $f_0$ follow.

Let $h_1(x)=f(x)-R_n(x)$ and $h_2(x)=x^{n+1}$.
Then by ({\ref{SerOff(x)}),
\begin{eqnarray}\label{SerOffn}
h_1(x)=\sum_{k=n+1}^{\infty}a_kx^k=x^{n+1}\sum_{k=0}^{\infty}a_{n+k+1}x^k, ~
f_n(x)=\frac{h_1(x)}{h_2(x)}=\sum_{k=0}^{\infty}a_{n+k+1}x^k,
\end{eqnarray}
$h_1^{(m)}(0^+)=h_2^{(m)}(0)=0$ for
$m\in\IN\cup\{0\}$ with $0\leq m \leq n$, and by differentiation,
$$\frac{h_1^{(n+1)}(x)}{h_2^{(n+1)}(x)}=\frac{f^{(n+1)}(x)}{(n+1)!},$$
so that$f_n$ has the same monotonicity
property as that of $f^{(n+1)}$ by \citep[Theorem 1.25]{2}. Hence the monotonicity properties
of $f_{2n}$ and $f_{2n-1}$ follow from Theorem \ref{th1}(2).

Next, differentiation gives
\begin{eqnarray}\label{fn'}
f_n'(x)&=&\frac{1}{x^{n+2}}\left\{x\left[f'(x)-\sum_{k=0}^nka_kx^{k-1}\right]-
(n+1)\left[f(x)-\sum_{k=0}^na_kx^k\right]\right\}\nonumber\\
&=&\frac{1}{x^{n+2}}\left[x
f'(x)-(n+1)f(x)+\sum_{k=0}^n(n+1-k)a_kx^k\right]=\frac{h_3(x)}{h_4(x)},
\end{eqnarray}
where $h_4(x)=x^{n+2}$ and
\begin{eqnarray*}
h_3(x)&=&xf'(x)-(n+1)f(x)+\sum_{k=0}^n(n+1-k)a_kx^k\\
&=&\sum_{k=n+2}^{\infty}(k-n-1)a_kx^k=x^{n+2}\sum_{k=0}^{\infty}(k+1)a_{n+k+2}x^k
\end{eqnarray*}
by (\ref{SerOff(x)}). It is easy to verify that $h_3^{(n+1)}(x)=xf^{(n+2)}(x)$. Hence $h_3^{(m)}(0^+)=h_4^{(m)}(0)=0$ for
$m\in\IN\cup\{0\}$ with $0\leq m \leq n+1$, and
\begin{eqnarray}\label{(n+2)-Df}
\frac{h_3^{(n+1)}(x)}{h_4^{(n+1)}(x)}=\frac{f^{(n+2)}(x)}{(n+2)!},
\end{eqnarray}
which shows that $f_n'$ has the same monotonicity
property as that of $f^{(n+2)}$ by \citep[Theorem 1.25]{2} and (\ref{fn'}). Consequently, the convexity (
concavity ) of $f_{2n}$ ( $f_{2n+1}$, respectively ) follows from
Theorem \ref{th1}(2).

The limiting value $f_n(0^+)=a_{n+1}$ follows from (\ref{SerOffn}). By the definition of $f_n$,
$f_n(1/2)=A_n$. Hence the double inequality (\ref{Ineq2OfR}) and its equality case follow from the
monotonicity and concavity properties of
$f_{2n-1}$ and (\ref{(cn+1)+(An+1)}). Taking $n=1$ in (\ref{Ineq2OfR}), we obtain (\ref{Ineq3OfR})
and its equality case.

(2) Let $h_5(x)=f(x)-S_n(x)=(1-2x)^{2(n+1)}\sum_{k=0}^{\infty}b_{n+k+1}(1-2x)^{2k}$ and
$h_6(x)=(1-2x)^{2(n+1)}$. Then $h_5^{(m)}(1/2)=h_6^{(m)}(1/2)=0$ for
$m\in\IN\cup\{0\}$ with $0\leq m \leq 2n+1$, and
\begin{eqnarray}\label{SerOfgn}
g_n(x)=\frac{h_5(x)}{h_6(x)}=\sum_{k=0}^{\infty}b_{n+k+1}(1-2x)^{2k}.
\end{eqnarray}
By the definition of $g_n$, $g_n(0^+)=h_5(0^+)=c_n$, and by (\ref{SerOfgn}), $g_n((1/2)^-)=b_{n+1}$. Hence it follows from (\ref{SerOfgn}) and Lemma \ref{Lem1}(4) that $g_n$ is strictly increasing and concave from $(0,1/2)$ onto $(c_n,b_{n+1})$.

Next, applying the method used to prove (\ref{SerOf(2n-1)-Df}) and (\ref{SerOf(2n)-Df}), we can obtain the following derivatives
\begin{eqnarray}\label{(2m+1)-Dgn}
g_n^{(2m+1)}(x)&=&-2^{2m+1}\sum_{k=1}^{\infty}\frac{(2m+2k)!}{(2k-1)!}b_{n+m+k+1}(1-2x)^{2k-1},\label{(2m+1)-Dgn}\\
g_n^{(2m)}(x)&=&2^{2m}\sum_{k=0}^{\infty}\frac{(2m+2k)!}{(2k)!}b_{n+m+k+1}(1-2x)^{2k}.\label{(2m)-Dgn}
\end{eqnarray}
Hence by Lemma \ref{Lem1}(4), $g_n^{(2m)}$ ( $g_n^{(2m+1)}$ ) is strictly increasing (
decreasing, respectively ) on $(0,1/2)$ for $m\in\IN\cup\{0\}$. This also yields the concavity ( convexity ) property of $g_n^{(2m)}$ $\left(
g_n^{(2m+1)}, \mbox{respectively} \right)$.

The proof of the complete monotonicity property for $G_{n,\thinspace m}$ is similar to that for $F_{1,\,n}$ in Theorem \ref{th1}(2), and we omit the details. By (\ref{(2m+1)-Dgn}) and (\ref{(2m)-Dgn}), we obtain the following limiting values
\begin{eqnarray*}
G_{n,\,2m+1}(0^+)&=&-2^{2m+1}\sum_{k=1}^{\infty}\frac{(2m+2k)!}{(2k-1)!}b_{n+m+k+1}, ~G_{n,\,2m+1}((1/2)^-)=0,\\
G_{n,\,2m}((1/2)^-)&=&-2^{2m}\sum_{k=0}^{\infty}\frac{(2m+2k)!}{(2k)!}b_{n+m+k+1}, ~G_{n,\,2m+1}(0^+)=-4^m(2m)!b_{n+m+1}.
\end{eqnarray*}

The double inequality (\ref{Ineq4OfR}) and its equality case follow from the
monotonicity and concavity properties of $g_n$. Taking $n=1$ in (\ref{Ineq4OfR}), we obtain (\ref{Ineq5OfR}) and its equality case. $\Box$

\begin{corollary}\label{Col2}
Let $a_n$, $b_n$, $c_n$ and $A_n$ are as in (\ref{an})--(\ref{cn}). Then we have the following conclusions:

(1) For all $n\in\IN\cup\{0\}$, $c_n<b_{n+1}$.

(2) $A_{2n+1}$ ( $A_{2n}$ ) is strictly increasing (decreasing, respectively) in $n\in\IN$, and for all $n\in\IN$,
\begin{eqnarray}\label{Sign1OfAn}
a_{2n+2}<A_{2n+1}<\lim_{n\to\infty}A_n=0<a_{2n+1}+\frac12 a_{2n+2}< A_{2n}<a_{2n+1}.
\end{eqnarray}

(3) $|a_n|$ is strictly decreasing in $n\in\IN\setminus\{1\}$. In particular, for $n\in\IN\setminus\{1\}$ and $\mu=65/108=0.601851\cdots$,
\begin{eqnarray}
\lim_{n\to\infty}a_n=0<-a_{2n+2}&<&a_{2n+1}<-a_{2n}<\cdots<a_3<a_2<-a_1,\label{Ineq1Ofan}\\
\left(1-\frac{1-\mu}{2^{n-1}-\mu}\right)\zeta(n)&<&\zeta(n+1)<\frac{\zeta(n)+1}{2}<\zeta(n).\label{IneqOfzeta}
\end{eqnarray}
\end{corollary}

{\it Proof.} (1) The inequality in part (1) holds by the result for $g_n$ in Theorem \ref{th2}(2).

(2) It follows from the monotonicity properties of $f_{2n+1}$ and $f_{2n}$ stated in Theorem \ref{th2}(1) that
\begin{eqnarray}\label{Ineq1OfA2n}
A_{2n}<a_{2n+1} \mbox{ ~and ~} a_{2n+2}< A_{2n+1}, ~n\in\IN.
\end{eqnarray}
By (\ref{Id1}), $\lim_{n\to\infty}a_n=0$. From (\ref{(cn+1)+(An+1)}) and (\ref{Ineq1OfA2n}), we obtain
\begin{eqnarray}
a_{2n+2}&<& A_{2n+1}=2(A_{2n}-a_{2n+1})<0,\label{Sign2OfAn}\\
a_{2n+1}&>&A_{2n}=a_{2n+1}+\frac12 A_{2n+1}>a_{2n+1}+\frac12 a_{2n+2},\label{Sign3OfAn}
\end{eqnarray}
and hence $\lim_{n\to\infty}A_n=0$ by Pinching Theorem for limits.

It follows from (\ref{zeta}), (\ref{lam-beta}) and (\ref{SignOfan}) that
\begin{eqnarray}\label{Ineq2Ofan}
a_{2n+1}+\frac12 a_{2n+2}
=\zeta(2n+3)+\zeta(2n+1)-2\eta(2n+2)=\sum_{k=2}^{\infty}\frac{[k+(-1)^k]^2}{k^{2n+3}}>0.
\end{eqnarray}
Hence (\ref{Sign1OfAn}) holds by (\ref{Ineq1OfA2n})--(\ref{Ineq2Ofan}).

Next, by (\ref{(cn+1)+(An+1)}) and (\ref{Sign2OfAn}), we obtain
\begin{eqnarray*}
A_{2n+1}-A_{2n+3}&=&2[(A_{2n}-a_{2n+1})-(A_{2n+2}-a_{2n+3})]\\
&=&2[A_{2n}-2(A_{2n+1}-a_{2n+2})+a_{2n+3}-a_{2n+1}]\\
&=&2[2a_{2n+2}+a_{2n+3}-3(A_{2n}-a_{2n+1})]\\
&<& 2\left(2a_{2n+2}+a_{2n+3}-\frac32 a_{2n+2}\right)=a_{2n+2}+2 a_{2n+3},
\end{eqnarray*}
and hence by (\ref{SignOfan}), (\ref{zeta}) and by (\ref{lam-beta}),
\begin{eqnarray}\label{Ineq4OfAn}
A_{2n+1}-A_{2n+3}&<& 2[\zeta(2n+3)-\zeta(2n+1)]+4[\zeta(2n+3)-\eta(2n+4)]\nonumber\\
&=&2[3\zeta(2n+3)-\zeta(2n+1)-2\eta(2n+4)]\nonumber\\
&=&-2\sum_{k=3}^{\infty}\frac{k^3-3k-2(-1)^k}{k^{2n+4}}
<-8\sum_{k=3}^{\infty}\frac{k+1}{k^{2n+4}}<0
\end{eqnarray}
since $k^3-3k\geq 4k+6\geq4(k+1)+2(-1)^k$ for $k\geq3$. This yields the monotonicity of $A_{2n+1}$.

By (\ref{Ineq1OfA2n}) and (\ref{Sign3OfAn}), we have
$$A_{2n+2}-A_{2n}<a_{2n+3}-a_{2n+1}-\frac12 a_{2n+2},$$
and hence it follows from (\ref{zeta}), (\ref{lam-beta}) and (\ref{SignOfan}) that
\begin{eqnarray*}
A_{2n+2}-A_{2n}&<&2[\zeta(2n+3)-\eta(2n+4)-\zeta(2n+1)+\eta(2n+2)]-\zeta(2n+3)+\zeta(2n+1)\\
&=&\zeta(2n+3)-\zeta(2n+1)+2[\eta(2n+2)-\eta(2n+4)]\\
&=&-\sum_{k=2}^{\infty}\frac{k^2-1}{k^{2n+4}}\left[k+2(-1)^k\right]<0.
\end{eqnarray*}
This yields the the monotonicity of $A_{2n}$.

(3) The limiting value of $a_n$ follows from (\ref{Id1}). Computation gives: $a_4=-0.395066\cdots$, $a_3=0.510048\cdots$, $a_2=2.404113\cdots$ and $a_1=-2.644934\cdots$. Hence $-a_4<a_3<a_2<-a_1$.

By (\ref{SignOfan}), in order to prove the monotonicity of $|a_n|$, we need only to prove (\ref{Ineq1Ofan}). By applying the Monotone l'H\^opital's Rule \citep[Theorem 1.25]{2}, one can  easily verify that the function $P_1(x)\equiv [(x+1)^8-1]/[x(x+1)^7]$ is strictly decreasing from $(0,1/2]$ onto $[6305/2187,8)$, by which we see that
$$P_2(m)\equiv 2m+1-2m\left(\frac{2m}{2m+1}\right)^7=\frac{[1+(2m)^{-1}]^8-1}{(2m)^{-1}[1+(2m)^{-1}]^7}=P_1\left(\frac{1}{2m}\right)$$
is strictly increasing in $m\in\IN$ with $P_2(1)=\omega\equiv 6305/2187=2.882944\cdots$ and $P_2(\infty)=8$.
Hence it follows from (\ref{SignOfan}), (\ref{zeta})
and (\ref{lam-beta}) that for $n\in\IN\setminus\{1\}$,
\begin{eqnarray}\label{Ineq3Ofan}
\frac{a_{2n+1}+a_{2n+2}}{2}&=&\zeta(2n+3)-\eta(2n+2)=\sum_{k=2}^{\infty}\frac{(-1)^kk+1}{k^{2n+3}}\nonumber\\
&=&\frac{1}{2^{2n+3}}\left[3-2\left(\frac{2}{3}\right)^{2n+3}\right]+\frac{1}{4^{2n+3}}\left[5-4\left(\frac{4}{5}\right)^{2n+3}\right]
+\cdots\nonumber\\
&{}&+\frac{1}{(2m)^{2n+3}}\left[(2m+1)-2m\left(\frac{2m}{2m+1}\right)^{2n+3}\right]+\cdots\nonumber\\
&\geq&
\sum_{m=1}^{\infty}\frac{P_2(m)}{(2m)^{2n+3}}>\frac{\omega}{2^{2n+3}}\sum_{m=1}^{\infty}\frac{1}{m^{2n+3}}=\frac{\omega}{2^{2n+3}}\zeta(2n+3),
\end{eqnarray}
and
\begin{eqnarray}\label{Ineq4Ofan}
\frac{a_{2n+1}+a_{2n}}{2}&=&2\zeta(2n+1)-\eta(2n+2)-\zeta(2n-1)=\sum_{k=1}^{\infty}\frac{2k-k^3+(-1)^k}{k^{2n+2}}\nonumber\\
&<&-\sum_{k=2}^{\infty}k^{-2n-2}\left(\frac{k^2}{2}+\frac{k}{2}\right)
=1-\frac{\zeta(2n)+\zeta(2n+1)}{2}<0
\end{eqnarray}
since $k^3-k^2/2-5k/2\geq1\geq(-1)^k$ for $k\geq2$.
Hence it follows from (\ref{SignOfan}), (\ref{Ineq3Ofan}) and (\ref{Ineq4Ofan}) that $0<-a_{2n+2}<a_{2n+1}<-a_{2n}$,
so that (\ref{Ineq1Ofan}) holds.

Similarly to (\ref{Ineq3Ofan}), we can easily obtain
\begin{eqnarray*}
\zeta(n+1)-\eta(n)&=&\sum_{k=2}^{\infty}\frac{(-1)^kk+1}{k^{n+1}}>\frac{1}{2^{n+1}}\sum_{m=1}^{\infty}\frac{P_3(m)}{m^{n+1}}
\\
&>&\frac{P_3(1)}{2^{n+1}}\sum_{m=1}^{\infty}\frac{1}{m^{n+1}}=2^{1-n}\mu\zeta(n+1),
\end{eqnarray*}
for $n\in\IN$, where $P_3(m)=2m+1-2m[2m/(2m+1)]^3$ with $P_3(1)=65/27=4\mu$. Hence by (\ref{lambda(n+1)}),
$$\zeta(n+1)>\frac{\eta(n)}{1-2^{1-n}\mu}=\frac{1-2^{1-n}}{1-2^{1-n}\mu}\zeta(n)=\left(1-\frac{1-\mu}{2^{n-1}-\mu}\right)\zeta(n)$$
for $n\in\IN\setminus\{1\}$. This yields the first inequality in (\ref{IneqOfzeta}). On the other hand, we have
\begin{eqnarray*}
\zeta(n+1)-\zeta(n)=-\sum_{k=2}^{\infty}\frac{k-1}{k^{n+1}}<-\frac12\sum_{k=2}^{\infty}\frac{1}{k^n}
=\frac{1-\zeta(n)}{2}
\end{eqnarray*}
for $n\in\IN\setminus\{1\}$, which yields the second inequality in (\ref{IneqOfzeta}). The third inequality in (\ref{IneqOfzeta}) is clear. $\Box$

\section{\normalsize Proof of Theorem \ref{th3}}

(1) Clearly, $F_1(1/2)=(5\log2)/\pi$. By (\ref{SerOfR}) and (\ref{SerOfB}), it is easy to see that $F_1(0^+)=1$.

Let $f$ be as in Theorem \ref{th1}. Then $F_1(x)=[f(x)/B(x)]+1$. Differentiation gives
\begin{eqnarray}\label{F1'}
\frac{B(x)F_1'(x)}{4(1-2x)}=H(x)\equiv
\frac{f'(x)+f(x)H_1(x)}{4(1-2x)},
\end{eqnarray}
where $H_1(x)=\psi(1-x)-\psi(x)$. Clearly, $H_1(1/2)=0$ and
$H_1(0^+)=\infty$.

By \citep[6.3.7, 4.3.67 \&23.2.16]{1} and (\ref{lambda(n+1)}), we have
\begin{eqnarray}
H_1(x)&=&\psi(1-x)-\psi(x)=\pi \cot(\pi x)=\pi\tan\left(\pi\left(\frac12-x\right)\right) ,\label{psi(1-x)}\\
\tan z&=&2\sum_{n=1}^{\infty}\left(\frac{2}{\pi}\right)^{2n}\lambda(2n)z^{2n-1}, ~|z|<\frac{\pi}{2}.\label{tan}
\end{eqnarray}
It follows from (\ref{psi(1-x)}) and (\ref{tan}) that
\begin{eqnarray}\label{SerOfH1}
H_1(x)=\pi\tan\left(\pi\left(\frac12-x\right)\right)
=-\sum_{k=1}^{\infty}2^{2k+1}\lambda(2k)\left(x-\frac12\right)^{2k-1}=4\sum_{k=1}^{\infty}\lambda(2k)(1-2x)^{2k-1}
\end{eqnarray}
by which we obtain the following limiting values
\begin{eqnarray}\label{n-DH1(1/2)}
H_1^{(2n+1)}(1/2)=-2^{2n+3}(2n+1)!\lambda(2n+2) ~\mbox{and }
~H_1^{(2n)}(1/2)=0, ~n\in\IN\cup\{0\}.
\end{eqnarray}

Let $H_2(x)=f(x)H_1(x)$. Then by differentiation and the Leibniz formula,
\begin{eqnarray}\label{n-DH2}
H_2^{(n)}(x)=\sum_{k=0}^nC_n^kf^{(k)}(x)H_1^{(n-k)}(x),
\end{eqnarray}
where $C_n^k=n!/[(n-k)!k!]$. By Theorem \ref{th1}(1),
$f^{(2n+1)}(1/2)=0$ and $f^{(2n)}(1/2)=4^n(2n)!b_n$ for $n\in\IN\cup\{0\}$. Hence
from (\ref{n-DH1(1/2)}) and (\ref{n-DH2}), we obtain the following values
\begin{eqnarray*}
H_2^{(2n)}\left(\frac12\right)&=&\sum_{k=0}^{2n}C_{2n}^kf^{(k)}\left(\frac12\right)H_1^{(2n-k)}\left(\frac12\right)\\
&=&\sum_{m=1}
^nC_{2n}^{2m-1}f^{(2m-1)}\left(\frac12\right)H_1^{(2n-2m+1)}\left(\frac12\right)+\sum_{m=0}^nC_{2n}^{2m}f^{(2m)}\left(\frac12\right)H_1^{(2n-2m)}\left(\frac12\right)=0,\\
H_2^{(2n+1)}\left(\frac12\right)&=&\sum_{k=0}^{2n+1}C_{2n+1}^kf^{(k)}\left(\frac12\right)H_1^{(2n-k+1)}\left(\frac12\right)\\
&=&\sum_{m=1}^{n+1}C_{2n+1}^{2m-1}f^{(2m-1)}\left(\frac12\right)H_1^{(2n-2m+2)}\left(\frac12\right)+\sum_{m=0}^nC_{2n+1}^{2m}f^{(2m)}\left(\frac12\right)H_1^{(2n-2m+1)}\left(\frac12\right)\\
&=&\sum_{m=0
}^nC_{2n+1}^{2m}f^{(2m)}\left(\frac12\right)H_1^{(2n-2m+1)}\left(\frac12\right)=-(2n+1)!2^{2n+3}\sum_{m=0 }^nb_m\lambda(2n-2m+2)
\end{eqnarray*}
for $n\in\IN\cup\{0\}$. Therefore, $H_2(x)$ has the following power series expansion
\begin{eqnarray}\label{SerOfH2}
H_2(x)=\sum_{n=0}^{\infty}\frac{H_2^{(2n+1)}(1/2)}{(2n+1)!}\left(x-\frac{1}{2}\right)^{2n+1}
=4\sum_{n=0}^{\infty}\left[\sum_{k=0
}^nb_k\lambda(2n-2k+2)\right](1-2x)^{2n+1}.
\end{eqnarray}

Let $d_n$ be as in Corollary \ref{Col1}. Then it follows from (\ref{SerOff(x)}), (\ref{F1'}) and (\ref{SerOfH2}) that
\begin{eqnarray}
H(x)&=&\frac{4\sum_{n=0}^{\infty}\left[\sum_{k=0
}^nb_k\lambda(2n-2k+2)\right](1-2x)^{2n+1}-4\sum_{n=1}^{\infty}nb_n(1-2x)^{2n-1}}{4(1-2x)}\nonumber\\
&=&\sum_{n=0}^{\infty}\left[\sum_{k=0
}^nb_k\lambda(2n-2k+2)-(n+1)b_{n+1}\right](1-2x)^{2n}
=\sum_{n=0}^{\infty}d_n(1-2x)^{2n}.\label{SerOfH}
\end{eqnarray}

Clearly, $H((1/2)^-)=d_0=-0.290171\cdots$, and
$H(0^+)=[a_1+a_0H_1(0^+)]/4=\infty$ by (\ref{F1'}) and Theorem \ref{th1}. It follows from Corollary \ref{Col1}(3) and (\ref{SerOfH}) that $H$ is strictly
decreasing and convex from $(0,1/2)$ onto $(d_0,\infty)$. Hance $H$ has a unique zero $x_1\in(0,1/2)$ such that $H(x)>0$ for $0<x<x_1$, and $H(x)<0$ for $x_1<x<1/2$, so that the
piecewise monotonicity property of $F_1$ follows from (\ref{F1'}).

By (\ref{B}), (\ref{F1'}) and (\ref{psi(1-x)}), we have
\begin{eqnarray*}
4(1-2x)H(x)&=&(1-2x)R(x)+[1+x(1-x)]R'(x)+\frac{\pi^2\cos(\pi x)}{\sin^2(\pi
x)}\\
&{}&+\left\{[1+x(1-x)]R(x)-\frac{\pi}{\sin(\pi x)}\right\}H_1(x),
\end{eqnarray*}
and hence
$$H\left(\frac14\right)=\frac12\left\{3\log2+\frac{19}{16}R'\left(\frac14\right)+\pi^2\sqrt2
+\pi\left[\frac{19}{16}R\left(\frac14\right)-\pi\sqrt2\right]\right\}=0.095698\cdots$$
 by (\ref{R(1/4)}) and (\ref{R'(1/4)}). This shows that $x_1\in(1/4,1/2)$.

Next, it follows from (\ref{psi(1-x)}) that
\begin{eqnarray}\label{B'}
B'(x)=-B(x)H_1(x) \mbox{ ~and ~} H_1'(x)=-\frac{\pi^2}{\sin^2(\pi x)}=-B(x)^2.
\end{eqnarray}
By (\ref{F1'}) and (\ref{SerOfH}),
$$F_1'(x)=\frac{4(1-2x)H(x)}{B(x)}=\frac{4}{B(x)}\sum_{n=0}^{\infty}d_n(1-2x)^{2n+1},$$
and hence by (\ref{B'}) and differentiation,
\begin{eqnarray}\label{F1''}
F_1''(x)=\frac{4}{B(x)}\left[H_1(x)\sum_{n=0}^{\infty}d_n(1-2x)^{2n+1}-2\sum_{n=0}^{\infty}(2n+1)d_n(1-2x)^{2n}\right].
\end{eqnarray}
It follows from (\ref{n-DH1(1/2)}), (\ref{F1''}) and Corollary \ref{Col1}(3) that
\begin{eqnarray}\label{F1''(1/2)}
F_1''(1/2)=-8d_0/\pi >0.
\end{eqnarray}
On the other hand, by (\ref{F1'}) and (\ref{psi(1-x)}),
\begin{eqnarray*}
F_1'(x)=\frac{f'(x)+f(x)H_1(x)}{B(x)}=\frac{f'(x)\sin(\pi x)+\pi f(x)\cos(\pi x)}{\pi},
\end{eqnarray*}
and by differentiation and (\ref{psi(1-x)}), we obtain
\begin{eqnarray}\label{1F1''}
F_1''(x)=\frac{1}{\pi}\left[f''(x)\sin(\pi x)+2\pi f'(x)\cos(\pi x)-\pi^2 f(x)\sin(\pi x)\right],
\end{eqnarray}
and hence we have
\begin{eqnarray}\label{F1''(0)}
F_1''(0^+)=2f'(0^+)=2a_1=-(2+\pi^2/3)<0.
\end{eqnarray}
By (\ref{F1''(1/2)}) and (\ref{F1''(0)}), we see that $F_1$ is neither convex nor concave on $(0,1/2)$.

The double inequality (\ref{Ineq6OfR}) is clear. Since $F_1(0^+)=1$, the coefficient $\alpha=1$ is best
possible.

It follows from (\ref{R(1/4)}) and the piecewise monotonicity property of $F_1$ that
$$\delta=F_1(x_1)>F_1\left(\frac14\right)=\frac{19R(1/4)}{16\pi\sqrt2}=\frac{19\log8}{8\pi\sqrt2 }=1.111592\cdots,$$
and hence the first inequality in (\ref{delta}) holds. On the other hand,
it follows from (\ref{Ineq5OfR}) and (\ref{B}) that
\begin{eqnarray}\label{Ineq11OfR}
R(x)\leq \frac{B(x)}{1+x(1-x)}[1+h_8(x)]
\end{eqnarray}
for all $x\in(0,1/2]$, where $h_8(x)=\frac{1}{\pi}S_2(x)\sin(\pi x)$ and $S_2(x)=b_0+b_1(1-2x)^2+b_2(1-2x)^4$. Let $h_{10}(x)=b_1+2b_2(1-2x)^2$ and $h_{11}(x)=(1-2x)\tan(\pi x)$. Then by differentiation,
\begin{eqnarray}\label{h8'}
h_8'(x)/\cos(\pi x)=h_9(x)\equiv
S_2(x)-4h_{10}(x)h_{11}(x)/\pi.
\end{eqnarray}
Clearly, $h_9(0)=b_0+b_1+b_2>0$ and $h_9((1/2)^-)=b_0-8b_1/\pi^2=-0.235204\cdots$. Since $b_2<0$ and since
$b_1+2b_2=0.662442\cdots$, $h_{10}$ is
strictly increasing from $[0,1/2]$ onto $[b_1+2b_2,b_1]$ and
$$S_2'(x)/[4(1-2x)]=-h_{10}(x)<-(b_1+2b_2)<0, ~0<x\leq 1/2,$$
so that $S_2$ is strictly decreasing from $[0,1/2]$ onto
$[b_0,b_0+b_1+b_2]$. It is easy to show that the function
$h_{11}$ is strictly increasing from $(0,1/2)$
onto $(0,2/\pi)$. Hence by (\ref{h8'}), $h_9$ is strictly decreasing from
$(0,1/2)$ onto $(b_0-8b_1/\pi^2, b_0+b_1+b_2)$, and $h_9$ has a unique zero $x_0\in(0,1/2)$ such that $h_8'(x)>0$ for $x\in(0,x_0)$, and $h_8'(x)<0$ for $x\in(x_0,1/2]$. Computation gives:
\begin{eqnarray*}
h_9(0.276937)=0.00000007895\cdots, ~h_9(0.276938)=-0.00000137425\cdots
\end{eqnarray*}
showing that $x_0\in(0.276937,0.276938)$. Therefore
$$h_8(x)\leq h_8(x_0)<\frac{1}{\pi}S_2(0.276937)\sin(0.276938\pi)=0.11214596\cdots<0.112146,$$
which yields the second inequality in (\ref{delta}) by (\ref{Ineq11OfR}).

(2) Clearly, $F_2(x)=f(x)/\{[1+x(1-x)]B(x)\}$. By differentiation,
\begin{eqnarray}\label{F2'}
\frac{[1+x(1-x)]^2}{1-2x}B(x)F_2'(x)=H_3(x)\equiv
4[1+x(1-x)]H(x)-f(x),
\end{eqnarray}
where $H$ is as in (\ref{F1'}). It follows from (\ref{SerOff(x)}) and (\ref{SerOfH}) that
\begin{eqnarray}
H_3(x)&=&\left[5-(1-2x)^2\right]H(x)-f(x)\label{H31}\\
&=&\sum_{n=0}^{\infty}5d_n(1-2x)^{2n}-\sum_{n=0}^{\infty}d_n(1-2x)^{2(n+1)}-\sum_{n=0}^{\infty}b_n(1-2x)^{2n}\nonumber\\
&=&\sum_{n=0}^{\infty}(5d_n-b_n)(1-2x)^{2n}-\sum_{n=1}^{\infty}d_{n-1}(1-2x)^{2n}
=\sum_{n=0}^{\infty}D_n(1-2x)^{2n},\label{H3}
\end{eqnarray}
where $D_n$ is as in Corollary \ref{Col1}. Clearly, $H_3(0^+)=4H(0^+)-a_0=\infty$ and
$H_3(1/2)=D_0<0$. Hence by (\ref{H3}) and Corollary \ref{Col1}(4), it is clear that $H_3$ is strictly
completely monotonic from $(0,1/2)$ onto $(D_0,\infty)$. This shows that $H_3$ has a unique zero $x_2\in(0,1/2)$
such that $H_3(x)>0$ for $x\in(0,x_2)$, and
$H_3(x)<0$ for $x\in(x_2,1/2]$, and hence the piecewise monotonicity
property of $F_2$ follows from (\ref{F2'}).

Next, it follows from (\ref{B}), (\ref{F1'}), (\ref{F2'}) and (\ref{H3}) that
\begin{eqnarray}\label{2F2'}
F_2'(x)&=&\frac{[1+x(1-x)][f'(x)+f(x)H_1(x)]-(1-2x)f(x)}{[1+x(1-x)]^2B(x)}\nonumber\\
&=&\frac{[1+x(1-x)][f'(x)\sin(\pi x)+\pi f(x)\cos(\pi x)]-(1-2x)f(x)\sin(\pi x)}{\pi[1+x(1-x)]^2}\nonumber\\
&=&\frac{\sum_{n=0}^{\infty}D_n(1-2x)^{2n+1}}{[1+x(1-x)]^2B(x)}.
\end{eqnarray}
By (\ref{B'}) and the third equality in (\ref{2F2'}), and by differentiation,
\begin{eqnarray*}
F_2''(x)&=&-\frac{1}{[1+x(1-x)]^2B(x)}\left\{2\sum_{n=0}^{\infty}(2n+1)D_n(1-2x)^{2n}\right.\\
&{}&\left.+\left[2\frac{1-2x}{1+x(1-x)}
-H_1(x)\right]\sum_{n=0}^{\infty}D_n(1-2x)^{2n+1}\right\}
\end{eqnarray*}
from which we obtain
\begin{eqnarray}\label{F2''(1/2)}
F_2''(1/2)=-32D_0/(25\pi)=0.723202\cdots>0.
\end{eqnarray}

On the other hand, by the second equality in (\ref{2F2'}) and by differentiation,
\begin{eqnarray*}
\pi F_2''(x)&=&\frac{f''(x)\sin(\pi x)+2\pi f'(x)\cos(\pi x)-\pi^2 f(x)\sin(\pi x)}{1+x(1-x)}+\frac{2f(x)\sin(\pi x)}{[1+x(1-x)]^2}\\
&{}&-\frac{2(1-2x)}{[1+x(1-x)]^2}\left\{f'(x)\sin(\pi x)+\pi f(x)\cos(\pi x)
-\frac{(1-2x)f(x)\sin(\pi x)}{1+x(1-x)}\right\},
\end{eqnarray*}
by which we have
\begin{eqnarray}\label{F2''(0)}
F_2''(0^+)=2[f'(0^+)-f(0^+)]=2(a_1-a_0)=-(4+\pi^2/3)<0.
\end{eqnarray}

The assertion on the convexity and concavity property of $F_2$ now follows from (\ref{F2''(1/2)}) and (\ref{F2''(0)}).

(3) Since $F_3(x)=f(x)/[1+x(1-x)]$, the monotonicity property of $F_3$
follows from Theorem \ref{th1}(2).

Clearly, $F_3(1/2)=4f(1/2)/5=\log16-4\pi/5=\rho$, and $F_3(0^+)=f(0^+)=a_0=1$. For $x\in(0,1/2]$, let
$$F_4(x)=[1+x(1-x)]^2f''(x)-2(1-2x)[1+x(1-x)]f'(x)+2[2-3x(1-x)]f(x).$$
Clearly, the function $x\mapsto2-3x(1-x)$ is strictly decreasing from $[0,1/2]$ onto $[5/4,2]$. By Theorem \ref{th1}(2), $f(x)>0$, $f'(x)<0$ and $f''(x)>0$ for all
$x\in(0,1/2]$. Hence $F_4(x)>0$ for $x\in(0,1/2]$. Differentiation gives
\begin{eqnarray*}
[1+x(1-x)]^3F_3''(x)=F_4(x)>0
\end{eqnarray*}
for all $x\in(0,1/2]$, which yields the
convexity of $F_3$.

The double inequality (\ref{Ineq7OfR}) and its equality case are clear. $\Box$

\begin{remark}\label{Rem2}
We present the comparison of the bounds of $R(x)$ given in Theorems \ref{th1}--\ref{th3} as follows.

(1) The double inequality (\ref{Ineq1OfR}) is better than (\ref{Ineq7OfR}). In fact, it is clear that
$$\frac{b_0+(1-2x)P(x)}{1+x(1-x)}\geq\frac{b_0}{1+x(1-x)}\geq \frac{4}{5}b_0=\rho$$
for $x\in(0,1/2]$, with equality if and only if $x=1/2$. On the other hand, we have
\begin{eqnarray*}
\frac{b_0+(1-2x)Q(x)+B(x)}{1+x(1-x)}\leq\frac{b_0+(1-b_0)(1-2x)+B(x)}{1+x(1-x)}, ~x\in(0,1/2],
\end{eqnarray*}
with equality if and only if $x=1/2$.
Put $y=1-2x$. Then $y\in[0,1)$, $1+x(1-x)=(5-y^2)/4$, and
\begin{eqnarray*}
\frac{b_0+B(x)+(1-b_0)(1-2x)}{1+x(1-x)}\leq \rho+(1-\rho)(1-2x)+\frac{B(x)}{1+x(1-x)}\Leftrightarrow\\
4[b_0+(1-b_0)y]\leq \left(5-y^2\right)[\rho+(1-\rho)y]\Leftrightarrow (1-\rho)y^2+\rho y-1\leq0,
\end{eqnarray*}
which is true since the function $y\mapsto (1-\rho)y^2+\rho y-1$ is strictly increasing from $[0,1]$ onto $[-1,0]$.

(2) It is clear that the lower bound of $R(x)$ given in (\ref{Ineq1OfR}) is better than that given in (\ref{Ineq6OfR}).

(3) Let $F_5(x)=[b_0+(1-b_0)(1-2x)]\sin(\pi x)$ for $x\in(0,1/2]$. Then $F_5(0)=0$, $F_5(1/2)=b_0$, and
\begin{eqnarray}
b_0+B(x)+(1-b_0)(1-2x)\leq \delta B(x)\Leftrightarrow F_5(x)\leq \pi(\delta-1),\label{F51}\\
b_0+B(x)+(1-b_0)(1-2x)\geq \delta B(x)\Leftrightarrow F_5(x)\geq \pi(\delta-1).\label{F52}
\end{eqnarray}
Differentiation gives
\begin{eqnarray*}
\frac{F_5'(x)}{\cos(\pi x)}=F_6(x)\equiv \pi[b_0+(1-b_0)(1-2x)]-2(1-b_0)\tan(\pi x),
\end{eqnarray*}
which is clearly strictly decreasing from $(0,1/2)$ onto $(-\infty, \pi)$. Hence $F_6$ has a unique zero $x_5\in(0,1/2)$ such that $F_5$ is strictly increasing on $(0,x_5]$, and decreasing on $[x_5,1/2]$. By (\ref{delta}), $b_0<0.35057<\pi(\delta-1)<0.112146\pi=0.35231704\cdots<0.3523171$. Since
 $F_5(1/4)=(1+b_0)/(2\sqrt{2})=0.468155\cdots>0.112146\pi>\pi(\delta-1)$, there exist two numbers $x_6, x_7\in(0,1/2)$ with $x_6<x_7$ such that
\begin{eqnarray}\label{boundOfF5}
F_5(x)
\begin{cases}
<\pi(\delta-1), \mbox{ ~if ~} x\in(0,x_6)\cup(x_7,1/2],\\
=\pi(\delta-1), \mbox{ ~if ~} x=x_6 \mbox{ ~or ~} x_7,\\
>\pi(\delta-1), \mbox{ ~if ~} x\in(x_6, x_7).
\end{cases}
\end{eqnarray}
Consequently, by (\ref{F51})--(\ref{boundOfF5}), the first upper bound of $R(x)$ given in (\ref{Ineq1OfR}) is better (worse) than that given in (\ref{Ineq6OfR}) for $x\in(0,x_6]\cup[x_7,1/2]$ (for $x\in[x_6,x_7]$, respectively).

Next, let $F_7(x)=\left[b_0+b_1(1-2x)^2\right]\sin(\pi x)$. Then $F_7(0)=0$, $F_7(1/2)=b_0$, and
\begin{eqnarray}
b_0+b_1(1-2x)^2+B(x)\leq \delta B(x)\Leftrightarrow F_7(x)\leq \pi(\delta-1),\label{F71}\\
b_0+b_1(1-2x)^2+B(x)\geq \delta B(x)\Leftrightarrow F_7(x)\geq \pi(\delta-1).\label{F72}
\end{eqnarray}
By differentiation,
\begin{eqnarray*}
\frac{F_7'(x)}{\cos(\pi x)}=F_8(x)\equiv \pi\left[b_0+b_1(1-2x)^2\right]-4b_1h_{11}(x),
\end{eqnarray*}
where $h_{11}$ is as in (\ref{h8'}). Clearly, $F_8$ is strictly decreasing on $(0,1/2)$ with $F_8(0)=\pi(b_0+b_1)>0$ and $F_8((1/2)^-)=\pi b_0-8b_1/\pi=-0.738915\cdots$. Hence $F_8$ has a unique zero $x_8\in(0,1/2)$ such that $F_7$ is strictly increasing on $(0,x_8]$, and decreasing on $[x_8,1/2]$. Since
 $F_7(0.28)=0.352694\cdots>0.112146\pi>\pi(\delta-1)$, there exist two numbers $x_9, x_{10}\in(0,1/2)$ with $x_9<x_{10}$ such that
\begin{eqnarray}\label{boundOfF7}
F_7(x)
\begin{cases}
<\pi(\delta-1), \mbox{ ~if ~} x\in(0,x_9)\cup(x_{10},1/2],\\
=\pi(\delta-1), \mbox{ ~if ~} x=x_9 \mbox{ ~or ~} x_{10},\\
>\pi(\delta-1), \mbox{ ~if ~} x\in(x_9, x_{10}).
\end{cases}
\end{eqnarray}
Hence by (\ref{F71})--(\ref{boundOfF7}), the second upper bound of $R(x)$ given in (\ref{Ineq1OfR}) is better (worse) than that given in (\ref{Ineq6OfR}) for $x\in(0,x_9]\cup[x_{10},1/2]$ (for $x\in[x_9,x_{10}]$, respectively).

(4) Let $R_n$ and $S_n$ be as in Theorem \ref{th2}, and $\delta$ as in Theorem \ref{th3}. Then it follows from (\ref{Ineq1OfR}), (\ref{Ineq2OfR}), (\ref{Ineq4OfR}) and (\ref{Ineq6OfR}) that for all $n\in\IN$ and $x\in(0,1/2]$,
\begin{eqnarray}\label{Ineq10OfR}
\frac{B(x)+D_1(x)}{1+x(1-x)}\leq R(x)\leq \frac{B(x)+D_2(x)}{1+x(1-x)},
\end{eqnarray}
with equality in each instance if and only if $x=1/2$, where
\begin{eqnarray*}
D_1(x)&=&\max\left\{b_0+(1-2x)P(x), R_{2n+2}(x)+A_{2n+2}x^{2n+3}, S_{n+1}(x)+c_{n+1}(1-2x)^{2n+3}\right\},\\
D_2(x)&=&\min\left\{b_0+(1-2x)Q(x), R_{2n+1}(x)+A_{2n+1}x^{2(n+1)}, S_{n+1}(x), (\delta-1)
B(x)\right\}.
\end{eqnarray*}

(5) Applying \citep[Theorem 1.25]{2} and Theorem \ref{th3}(3), one can easily show that the function $F_9(x)\equiv[F_3(x)-1]/x$ is strictly increasing from $(0,1/2]$ onto $(-2-\pi^2/6,-2(1-\rho)]$, while $F_{10}(x)\equiv[F_3(x)-\rho]/(1-2x)$ is strictly decreasing from $(0,1/2)$ onto $(0,1-\rho)$, where $F_3$ is as in Theorem \ref{th3}.
\end{remark}

\begin{remark}\label{Rem3}
Theorems \ref{th1}--\ref{th3} show that the function $B(x)/[1+x(1-x)]=\pi/\{[1+x(1-x)]\sin(\pi x)\}$ is a very good approximation of $R(x)$.
\end{remark}

\begin{conjecture}\label{Conj1}
Let $F_3$ be as in Theorem \ref{th3}. Our computation supports the following conjecture: $F_3$ is strictly completely monotonic on $(0,1/2]$.

If this conjecture is true, then we can obtain sharp lower and upper bounds for $R(x)$, which are expressed in terms of $B(x)/[1+x(1-x)]$ and polynomials.
\end{conjecture}

\bigskip
\bigskip

X.-Y. Ma's email: mxy@zstu.edu.cn

T.-R. Huang's e-mail: htiren@zstu.edu.cn.

\clearpage

\end{document}